\documentclass[final,oneside,onecolumn]{article}
\usepackage[body={6in,9.5in},top=0.9in,left=1in,dvips]{geometry}
\usepackage{amsmath}
\usepackage{float}
\usepackage{pifont}
\usepackage{graphicx}
\usepackage{amsfonts}
\usepackage{amssymb}
\usepackage{amsthm}
\usepackage{mathbbol}
\usepackage{setspace}
\usepackage{paralist}
\usepackage{floatflt}
\usepackage{alltt}
\usepackage{color}
\usepackage{algorithm}
\usepackage{algorithmic}
\usepackage{hyperref}
\definecolor{string}{rgb}{0.7,0.0,0.0}
\definecolor{comment}{rgb}{0.13,0.54,0.13}
\definecolor{keyword}{rgb}{0.0,0.0,1.0}
\usepackage[parfill]{parskip}
\usepackage{color}

\usepackage{color-edits}
\addauthor{sh}{blue}
\addauthor{boom}{magenta}

\newcommand{\R}{\mathbb{R}}
\newcommand{\E}{\mathbb{E}}
\newcommand{\Prob}{\mathbb{P}}

\newtheorem{sublemma}{Lemma}[section]
\newtheorem{subdef}{Definition}[section]
\newtheorem{assp}{Assumption}[section]
\newtheorem{thm}{Theorem}[section]

\begin{document}

\renewcommand{\arraystretch}{0.5}

\title{\bf \large A New Convergence Analysis of Two Stochastic Frank-Wolfe Algorithms}
\date{April 2025}
\author{Natthawut Boonsiriphatthanajaroen and Shane G. Henderson\\School of Operations Research and Information Engineering, Cornell University}
\maketitle

\thispagestyle{empty}
\noindent

\begin{abstract}
We study the convergence properties of the original and away-step Frank-Wolfe algorithms for linearly constrained stochastic optimization assuming the availability of unbiased objective function gradient estimates. The objective function is not restricted to a finite summation form, like in previous analyses tailored to machine-learning applications. To enable the use of concentration inequalities we assume either a uniform bound on the variance of gradient estimates or uniformly sub-Gaussian tails on gradient estimates. With one of these regularity assumptions along with sufficient sampling, we can ensure sufficiently accurate gradient estimates. We then use a Lyapunov argument to obtain the desired complexity bounds, relying on existing geometrical results for polytopes. 

\end{abstract}
%assume two different bounded var or subGaus tails
%can add some notation
\section{Introduction}

%\bibitem{fw_non_cvx} \bibitem{trust_supermartingale}\bibitem{markov_bremaud}\bibitem{polyak}\bibitem{foss}\bibitem{fw_analysis}
%\bibitem{stoch_ls}\bibitem{first_tr}\bibitem{cubic}\bibitem{matrix}\bibitem{1sfw}\bibitem{rob}

The Frank-Wolfe algorithm~\cite{frankwolfe1956}, also known as the conditional gradient method, is a widely used approach for constrained optimization that relies on a linear minimization oracle over the feasible set. Frank-Wolfe algorithms avoid projections onto the feasible set, making them particularly useful relative to projection methods on problems where projection is expensive. For instance, as highlighted in \cite{pokutta2022conditional}, projection can be more expensive than solving a linear minimization oracle in feasible sets such as the nuclear norm ball, $\ell_p$-balls, the Birkhoff polytope, and the matching polytope. Examples of problems where projection complexity is particularly high include matrix completion and network routing. Frank-Wolfe algorithms have been applied to various machine learning tasks, such as training neural networks, and to regularized optimization problems~\cite{freundetal2020fwsum,zhangetal2020onesfw}.

We investigate these algorithms with a view towards simulation optimization with linear constraints that define a bounded polytope, i.e., we consider the problem
\begin{equation}\label{main_problem}
\begin{split}
\min_x & f(x) = \E[F(x,\xi)] \\
& \text{s.t. } Ax \le b
\end{split}
\end{equation}
where the set $\mathcal{X} = \{x: Ax \le b\}$ is a bounded polytope in $\R^d$ with $d < \infty$. In simulation applications, $\xi$ represents a random sample drawn from an underlying distribution and $F(x, \xi)$ represents an objective function realization obtained from a single simulation replication at the solution $x$. We can generate independent realizations of $F(x, \xi)$ and we estimate $f(x)$ using sample averages. We assume we have access to unbiased gradient estimates. Such estimates can be obtained in many, but not all, simulation applications using infinitesimal perturbation analysis, which can be computed, e.g., through automatic differentiation of simulation replications \cite{fuGradChapter2015}. In some applications gradient estimators can also be obtained using the likelihood ratio method or extensions thereof \cite{gly90,pengfu2025like}. 
The Frank-Wolfe algorithm iteratively solves a linear program (LP) using the gradient of the objective function at the current solution as the cost vector, and then moves towards an optimal vertex. When the optimal solution lies on a face of the feasible region, the algorithm can exhibit zig-zagging behavior between multiple vertices, leading to slower convergence. To mitigate this issue, the away-step Frank-Wolfe algorithm introduces an additional step that moves in the direction opposite to the least improving search direction, reducing zig-zagging and improving convergence. Theoretical convergence results for deterministic optimization problems, such as those in \cite{beckshtern2017linearly} and \cite{lacostejulienJaggi2015global}, rely on geometric properties of the feasible region. For deterministic, strongly convex objective functions, the away-step Frank-Wolfe algorithm has been shown to achieve faster convergence rates than the standard Frank-Wolfe algorithm.

We study the convergence behavior of both the standard Frank-Wolfe and the away-step Frank-Wolfe algorithms under certain function and gradient properties. Our goal is to find a solution $x$ such that $f(x) - f^{*} \le \varepsilon$ for some pre-specified $\varepsilon > 0$, where $f^{*}$ is the optimal objective value. A key challenge in our setting is that our gradient estimates are inaccurate due to Monte Carlo sampling error. This error can be controlled by increasing sample sizes at the cost of greater computational complexity. We analyze this tradeoff. Our analyses extend previous work that assumes a diminishing gradient variance or error \cite{mouline2015conv, gasnikov2024zoo} as an optimal solution is approached. Specifically, we establish both iteration and sample size complexity bounds that depend on the error tolerance $\varepsilon$; see Table~\ref{tab:complexity}. Moreover, we also establish that the (random) number of iterations required to achieve an $\varepsilon$-optimal solution has tails that decay at an exponential rate, i.e., that the required number of iterations is light-tailed. In addition to function properties such as strong convexity and Lipschitz continuity of gradients, we assume unbiased gradient estimates with either bounded variance or a uniform sub-Gaussian property.

The remainder of this paper is organized as follows. We place our work in context in Section~\ref{sec:related}. Section~\ref{sec:algs} specifies the algorithms. We prove preparatory results in Section~\ref{sec:prep}. Section ~\ref{sec:iteration} develops an iteration complexity analysis for both the standard and away-step Frank-Wolfe methods. Section~\ref{sec:samplesize} then establishes sample size complexity results under varying assumptions on the tail behavior of gradient estimates.

\section{Related Work}
\label{sec:related}

Convergence analyses of the Frank-Wolfe algorithm for deterministic objective functions have been extensively studied \cite{beckshtern2017linearly,lacostejulienJaggi2015affine,lacostejulienJaggi2015global}. The standard Frank-Wolfe algorithm requires $\mathcal{O}(\varepsilon^{-1})$ iterations to achieve $\varepsilon$-optimal solutions, while the away-step Frank-Wolfe algorithm improves to $\mathcal{O}(\log (\varepsilon^{-1}))$ iterations. The convergence rate of the away-step variant depends on the geometry of the feasible region, characterized by a property known as the pyramidal width. The pyramidal width does not need to be explicitly computed to be able to implement the algorithms or establish their complexity properties, which is fortunate because it is challenging to compute.

The study of the stochastic setting \cite{hazanluo2016var,freundlu2021fw,freundetal2020fwsum}, where the objective function has the form~\eqref{main_problem}, focuses on finite-sum minimization problems, where the objective function is the average of $n$ loss functions over a fixed dataset size $n$. Unlike the problem~\eqref{main_problem}, where the objective function $f$ is not inherently separable, these works assume a decomposition into multiple components $f_i$. Additionally, \cite{freundlu2021fw} and \cite{freundetal2020fwsum} consider objective functions of the form $f_i(x_i^T w)$, where $w$ is the optimization variable and $x_i$ is a given data point.

These studies propose a step-size policy that depends only on the iteration number. The gradient estimate is composed of $n$ entries, one for each data point. At each iteration, the gradient estimate is updated by randomly selecting a data point and updating only the corresponding gradient-estimate entry, eliminating the need for sample size analysis in gradient estimation. The iteration complexity in this setting is $\mathcal{O}(\varepsilon^{-1})$, with additional dependence on the dataset size $n$. Other works applying the Frank-Wolfe method, such as \cite{reddietal2016fw} and \cite{zhangetal2020onesfw}, focus on techniques like sampling methods and variance reduction rather than analyzing the convergence properties of the algorithm itself. In contrast, we explore the convergence properties of simulation-optimization problems where both the objective function and its gradient are estimated through sample averaging over multiple replications.

The quality of gradient approximations has been analyzed without relying on a finite-sum representation of the objective function in \cite{mouline2015conv}. There, it is assumed that the stochastic error in the gradient estimates decreases over iterations, leading to convergence results that guarantee a small error with high probability. A similar assumption is made in \cite{gasnikov2024zoo}, where the variance of the gradient estimate is assumed to diminish over time due to noise reduction. Under this assumption, the standard Frank-Wolfe algorithm achieves an iteration complexity of $\mathcal{O}(\varepsilon^{-1})$ when the variance diminishes to zero. However, these studies do not examine the required sample size for gradient estimation, which is a central question if Frank-Wolfe algorithms are to be applied in simulation optimization.

Related sample-complexity results have been established for other optimization algorithms, and partly inspired the present paper.
Works such as \cite{blanchetscheinberg2019trust}, \cite{menickellyscheinberg2018stoch}, and \cite{paqscheinberg2020ls} investigate optimization methods that incorporate simulated gradients to achieve an $\varepsilon$-neighborhood of the optimal solution, often using a Lyapunov function framework. In particular, the analysis of the trust-region method in \cite{blanchetscheinberg2019trust} centers on ``good events," where function or gradient estimates are sufficiently accurate. They show that the algorithm can be designed so that these good events occur with high probability, ensuring continued function improvement even in the presence of low-probability ``bad events." Moreover, they establish the necessary sample size required to make good events likely, providing important insights into sample complexity in stochastic optimization.

We establish convergence analyses and rate guarantees for the standard Frank-Wolfe and away-step Frank-Wolfe algorithms for the problem~\eqref{main_problem}, which extends the formulation studied in \cite{freundlu2021fw} and \cite{freundetal2020fwsum}. Our approach is partly inspired by the Lyapunov function analyses \cite{blanchetscheinberg2019trust, paqscheinberg2020ls} for trust-region methods. We confine our attention to strongly convex functions and assume the availability of unbiased gradient estimates. To understand sample size requirements we impose an assumption of either bounded variance, or a uniform sub-Gaussian property, of the simulated gradients. Compared to \cite{freundlu2021fw,mouline2015conv, freundetal2020fwsum}, and \cite{gasnikov2024zoo}, our work examines a more general class of objective functions than those seen in machine learning applications, namely those relevant in simulation optimization applications. We provide convergence analyses for the required sample size in gradient estimation. Table~\ref{tab:complexity} summarizes the complexity results established in this paper. The iteration complexity results are established in Theorem~\ref{thm:iteration_complexity}. The sample-size complexity results are established in Theorems~\ref{sample_size_thm} and \ref{th:SampleSizeSubG}. The total complexity results are the product of the sample size and iteration complexities. Moreover, Theorem~\ref{thm:iteration_complexity} shows that the (random) number of iterations required to reach an $\varepsilon$-optimal solution has a finite moment generating function in a neighborhood of zero, and therefore possesses tails that decay at an exponential rate.

\begin{table}[htb]
\begin{center}
\begin{tabular}{l  c c c c} 
 {\bf Algorithm} & {\bf Sub-Gaussian?} & {\bf Sample Size} & {\bf Iteration} & {\bf Total}\\ 
 \\
 Standard & No & $\mathcal{O}\left(\frac{1}{\varepsilon^4}\right)$ & $\mathcal{O}\left(\frac{1}{\varepsilon^2}\right)$ & $\mathcal{O}\left(\frac{1}{\varepsilon^6}\right)$ \\ 
 \\
 \\
 Away-Step & No & $\mathcal{O}\left(\frac{1}{\varepsilon^2}\right)$ & $\mathcal{O}\left(\frac{1}{\varepsilon}\right)$ & $\mathcal{O}\left(\frac{1}{\varepsilon^3}\right)$ \\
 \\
 \\
 Standard & Yes & $\mathcal{O}\left(\frac{1}{\varepsilon^2}\log\frac{1}{\varepsilon}\right)$ & $\mathcal{O}\left(\frac{1}{\varepsilon^2}\right)$ & $\mathcal{O}\left(\frac{1}{\varepsilon^4}\log\frac{1}{\varepsilon}\right)$ \\
 \\
 \\
 Away-Step & Yes & $\mathcal{O}\left(\frac{1}{\varepsilon}\log\frac{1}{\varepsilon}\right)$ & $\mathcal{O}\left(\frac{1}{\varepsilon}\right)$ & $\mathcal{O}\left(\frac{1}{\varepsilon^2}\log\frac{1}{\varepsilon}\right)$ \\
 \\

\end{tabular}
\end{center}
\caption{Complexity results established in this paper. The column ``Sub-Gaussian" indicates whether this is assumed of the gradient estimates or not. If it is not assumed, then we assume bounded variance. \label{tab:complexity}}
\end{table}

Compared to algorithms like Projected Stochastic Gradient Descent (PGD), the Frank-Wolfe algorithm has a similar iteration complexity but incurs an additional burden in terms of sample size. In \cite{lacostejulien2012pgd}, which is an extension of \cite{nemirovskilan2009robust}, it is shown that PGD has iteration complexity $\mathcal{O}({\varepsilon^{-1}})$ assuming strong convexity and $\mathcal{O}({\varepsilon^{-2}})$ with regular convexity. These results are similar to those for Stochastic Gradient Descent (SGD) on unconstrained problems, as discussed in \cite{garrigosgower2013conv}, for example, with the analysis for PGD exploiting the fact that the projection step is a contraction. Both these analyses (of stochastic PGD and SGD) rely on assumptions of unbiased gradients and bounded gradient variance, but under reasonable conditions only require a single replication when obtaining gradient estimates. This is a significant advantage relative to our approach, because we need to control the gradient estimation accuracy through sample-size selection. This advantage comes at the cost of projection; if projection is expensive then one might prefer to use Frank-Wolfe algorithms as analyzed here.

\section{Algorithms}
\label{sec:algs}
In this section, we present the main algorithms. We require the following assumptions.

\begin{assp}
The gradient of the objective function $f$ is $L-$Lipschitz continuous, i.e. 
$$||\nabla f(x) - \nabla f(y)|| \le L||x-y|| \ \ \forall x,y \in \mathcal{X}.$$
\end{assp}

\begin{assp}
The objective function $f$ is $\mu-$strongly convex, i.e., for some $\mu > 0$, 
$$f(y) \ge f(x) + \nabla f(x)^T(y-x) + \frac{\mu}{2}||x-y||^2 \ \ \forall x,y \in \mathcal{X}.$$
\end{assp}

\begin{assp}
We have access to an oracle that, at any feasible $x$, provides independent and identically distributed (i.i.d) vectors $G_1(x), G_2(x), \ldots$ such that $\E[G_1(x)] = \nabla f(x)$. 
\label{unbiased_g_assp}
\end{assp}

We apply a Frank-Wolfe algorithm to solve \eqref{main_problem}. In iteration $k$, the algorithm solves a linear program over the feasible region $\mathcal{X}$, where the objective function coefficients are given by a gradient estimate $g_k$. This gradient estimate is the average of multiple gradient realizations. In each iteration, we maintain a representation of the current solution as a convex combination of the set $\mathcal{V}$ of vertices of the feasible region, i.e., in iteration $k$, the solution $x_k$ at that iteration is represented as
$$x_k = \sum_{u \in \mathcal{V}} \alpha^{(k)}_u u,$$
where $\alpha^{(k)}_u$ are non-negative and sum to 1. This representation is not unique, but it plays a central role in the Away-step Frank-Wolfe algorithm's definition and analysis. We define the \textbf{active vertex set} at iteration $k$, $U_k$, to be the set of all vertices with positive weights, i.e. $U_k = \{v \in \mathcal{V}| \alpha^{(k)}_v > 0\}$. Let $N_k = |U_k|$ be the number of active vertices in iteration $k$.

We start with a feasible solution, and in each iteration, $k$, we determine a search direction then take a step of size $\gamma_k$ in the search direction. The step size $\gamma_k$ is limited to maintain feasibility. In the standard Frank-Wolfe algorithm, the maximum step size, $\gamma^{max}_k$, is 1, representing a step that exactly reaches a vertex of the feasible region. In the Away-step Frank-Wolfe algorithm, the maximum step size is determined by the current solution's representation, as discussed below. The step size that is taken, if not constrained by the maximal step size, is given by minimizing an upper bound on the objective function. Indeed, by convexity and $L-$Lipschitz continuity of the gradient, in the search direction $d_k$,
$$f(x_k + \gamma_k d_k) \le f(x_k) + \gamma_k\nabla f(x_k)^Td_k + \gamma^2_k\frac{L}{2}\|d_k\|^2.$$

Minimizing this quadratic in $\gamma_k$ and accounting for the maximal step size, we get the step size
$$\gamma_k = \min\left\{\gamma^{max}_k, \frac{-g_k^Td_k}{L\|d_k\|^2}\right\}.$$

\subsection{The Standard and Away-Step Frank-Wolfe Algorithms}

In the standard Frank-Wolfe algorithm, Algorithm~\ref{normal_fw_alg}, search always proceeds in the direction of a vertex. The algorithm may alternately search in the direction of multiple vertices when approaching an optimal point near a face of the polytope. This results in zig-zagging behavior, which may be more expensive than moving directly towards the face. This observation motivates Algorithm~\ref{away_fw_alg}, the Away-step Frank-Wolfe algorithm, which has the additional ability to move away from a ``bad" vertex, i.e., the vertex that produces the worst estimated function improvement. 

\begin{algorithm}[hbt!]
\caption{Frank-Wolfe}
\begin{algorithmic}
	\STATE Initialize a feasible solution $x_0 \in \mathcal{X}$
	\FOR{$k = 0,1,2,...$}
		\STATE Compute an estimate, $g_k$, of the gradient $\nabla f(x_k)$
		\STATE Compute $s_k = \arg\min_{s \in \mathcal{X}} g_k^Ts$ by solving an LP to obtain a vertex solution
		\STATE $d_k = s_k - x_k$
            \STATE $\gamma_k = \min\{1, \frac{\varepsilon}{2LD^2}\}$
		\STATE $x_{k+1} = x_k + \gamma_k d_k$
	\ENDFOR
\end{algorithmic}
\label{normal_fw_alg}
\end{algorithm}
\noindent

In the $k^{th}$ iteration of the Away-step Frank-Wolfe algorithm with current iterate $x_k$, we first obtain a gradient estimate, $g_k$, at $x_k$ and solve the LP to get the usual search direction. We then find another direction, the ``away direction," through the current iterate, {\em away} from an active vertex. In doing so, we select the active vertex $v_k$ that maximizes the dot product of $x_k-v_k$ with the negative gradient estimate. Then, we choose the better of the usual search direction and the away direction. To find the away direction vertex, we iterate over all current active vertices (those in the set $U_k$). As long as this set is modest in size, this step requires much less computation than solving the LP for the standard Frank-Wolfe direction. We initialize the active set, $U_0$, to a singleton set consisting of a vertex solution of the LP $\min_{x \in \mathcal{X}} \mathbb{1}^Tx$.

\begin{algorithm}[hbt!]
\caption{Away-Step Frank-Wolfe}
\begin{algorithmic}
	\STATE Initialize a feasible solution $x_0 \in \mathcal{X}$ and initialize the active set $U_0$ and coefficients appropriately
	\FOR{$k = 0,1,2,...$}
		\STATE Compute an estimate, $g_k$, of the gradient $\nabla f(x_k)$
		\STATE Compute $s_k = \arg\min_{s \in \mathcal{X}} g_k^Ts$ by solving an LP to obtain a vertex solution
		\STATE Compute $v_k = \arg\max_{u\in U_k} g_k^Tu$ 
		\STATE $d^{FW}_k = s_k - x_k$ \,\, (FW direction)
		\STATE $d^{A}_k = x_k - v_k$  \,\, (Away direction)
		
		\IF{$-g_k^Td^{FW}_k \ge -g_k^Td^{A}_k$}
			\STATE $d_k = d^{FW}_k$ and $\gamma^{max}_k = 1$ \,\, (normal FW step)
		\ELSE
			\STATE $d_k = d^{A}_k$ and set $\gamma^{max}_k$ as described in the text\,\, (away step)
		\ENDIF
		
		%\STATE $\alpha_k = \arg\min f(x_k + \gamma d_k)$ 
		\STATE $\gamma_k =  \min\left\{\gamma^{max}_k, \frac{-g_k^Td_k}{L\|d_k\|^2}\right\}$ 
		\STATE $x_{k+1} = x_k + \gamma_k d_k$
		\STATE Update the vertex representation of $x_{k+1}$ and active set as described in the text
	\ENDFOR
\end{algorithmic}
\label{away_fw_alg}
\end{algorithm}
%\noindent

The update of the vertex representation and active set depends on whether a normal Frank-Wolfe step or an away step is taken, and also on whether we take the maximal step \cite{lacostejulienJaggi2015global}. For a normal Frank-Wolfe step, the direction vector is $d_k = s_k - x_k$ for some $s_k \in \mathcal{V}$, and thus the maximal step size is $\gamma_k^{max}=1$. For an away-step, the maximum step size is restricted by the representation of the current iterate. An away-step gives $x_{k+1} = x_k + \gamma_k (x_k - v_k)$. Thus, for the vertex $v_k$ we have the update $\alpha^{(k+1)}_{v_k} = (1 + \gamma_k)\alpha^{(k)}_{v_k} - \gamma_k$. This weight must remain non-negative, so
$$0 \le \alpha^{(k+1)}_{v_k} = (1 + \gamma^{max}_k)\alpha^{(k)}_{v_k} - \gamma^{max}_k,$$
implying that
$$\gamma^{max}_k \le \frac{\alpha^{(k)}_{v_k}}{1 - \alpha^{(k)}_{v_k}}.
$$

Summarizing the above discussion, the updates to the active set and weights are as follows.
\begin{enumerate}
	\item Normal FW Step: 
		\begin{itemize}
			\item $\gamma_k = \gamma^{max}_k = 1$: Set $U_{k+1}= \{s_k\}$
			\item $\gamma_k < \gamma^{max}_k = 1$: Set $U_{k+1} = U_{k} \cup \{s_k\}$
			\item Update the weights $\alpha^{(k+1)}_{s_k} = (1 - \gamma_k)\alpha^{(k)}_{s_k} + \gamma_k$
			\item Update the weights $\alpha^{(k+1)}_{v} = (1 - \gamma_k)\alpha^{(k)}_{v}$ for all $v \in U_{k} \setminus \{s_k\}$
		\end{itemize}
	\item Away Step:
		\begin{itemize}
			\item $\gamma_k = \gamma^{max}_k$: Set $U_{k+1} = U_{k} \setminus \{v_k\}$
			\item $\gamma_k < \gamma^{max}_k$: Set $U_{k+1} = U_{k}$
			\item Update the weights $\alpha^{(k+1)}_{v_k} = (1 + \gamma_k)\alpha^{(k)}_{v_k} - \gamma_k$
			\item Update the weight $\alpha^{(k+1)}_{v} = (1 + \gamma_k)\alpha^{(k)}_{v}$ for all $v \in U_{k} \setminus \{v_k\}.$
		\end{itemize}
\end{enumerate}

\section{Preparatory Results for the Convergence Analysis}
\label{sec:prep}

We seek an $\varepsilon$-optimal point, for some fixed $\varepsilon > 0$, i.e., we run the algorithm until the condition $f(x) - f^{*} \le \varepsilon$ is satisfied, where $f^*$ is the optimal objective function value. To that end, define the terminal time $T_{\varepsilon}$ to be $T_{\varepsilon} = \inf\{k \ge 0 : f(x_k) - f^{*} \le \varepsilon\}$. We suppose, throughout this section, that $k < T_\varepsilon$, so we are yet to reach an $\varepsilon$-optimal point.

We begin with preparatory results assuming deterministic dynamics, which we leverage in a later section to obtain our main results. We adapt analyses from \cite{beckshtern2017linearly} and \cite{lacostejulienJaggi2015global} that assume exact gradients to our setting where there is an error in the gradient estimate. Since we assume deterministic dynamics in this section, we refer to a gradient {\em approximation} rather than a gradient {\em estimate}.

\subsection{Results for the Standard Frank-Wolfe Algorithm}

\begin{subdef}
Let $D$ be the diameter of the feasible region $\mathcal{X}$. We say that the vector $g_k$ is a {\em good gradient approximation} at iteration $k$ if
\[
\|g_k - \nabla f(x_k)\| \le \frac{\varepsilon}{4D}.
\]
\label{good_grad_standard}
\end{subdef}

\begin{sublemma}
If we have a good gradient approximation at iteration $k$, then
$$f(x_{k+1}) - f(x_k) \le -\beta_1\varepsilon,$$
where $\beta_1 = \min\{\frac{\varepsilon}{8LD^2}, \frac{1}{4}\}$.
\label{standard_fw_improvement}
\end{sublemma}

\begin{proof}
By convexity and Definition~\ref{good_grad_standard},
\begin{eqnarray*}
f(y) & \ge & f(x_k) + \nabla f(x_k)^T(y - x_k) \\
& \ge & f(x_k) + g^T_k(y - x_k) - \frac{\varepsilon}{4D}||y - x_k|| \\
& \ge & f(x_k) + g^T_k(y - x_k) - \frac{\varepsilon}{4}.
\end{eqnarray*}
Minimizing both sides, we have
$$f^{*} \ge f(x_k) + g^T_k(s_k - x_k) - \frac{\varepsilon}{4}.$$
Then, by convexity and the Lipschitz property,
\begin{eqnarray*}
f(x_{k+1}) & \le & f(x_k) + \gamma_k\nabla f(x_k)^T(s_k - x_k) + \frac{L\gamma^2_k}{2}||s_k-x_k||^2 \\
& \le & f(x_k) + \gamma_kg_k^T(s_k - x_k) + \gamma_k\frac{\varepsilon||s_k - x_k||}{4D} + \frac{L\gamma^2_k}{2}||s_k-x_k||^2 \\
& \le & f(x_k) + \gamma_k(f^{*} - f(x_k) + \frac{\varepsilon}{4}) + \gamma_k\frac{\varepsilon}{4} + \frac{L\gamma^2_k}{2}D^2. 
\end{eqnarray*}

Hence, we have the improvement
\begin{eqnarray*}
f(x_{k+1}) - f^{*} \le (1 - \gamma_k)(f(x_{k}) - f^{*}) + \frac{\gamma_k\varepsilon}{2} + \frac{L\gamma^2_kD^2}{2}.
\end{eqnarray*}

The step size policy $\gamma_k = \min\{\varepsilon/(2LD^2), 1\}$ corresponds to moving either to a point between the current solution and the vertex $s_k$ or moving all the way to the vertex $s_k$.

\textit{Case I}: Suppose $\gamma_k = \varepsilon / (2LD^2)$. Then, using the fact that the algorithm is not yet terminated, i.e. $f(x_k) - f^{*} > \varepsilon$,
\begin{eqnarray*}
f(x_{k+1}) - f^{*} & \le & (1 - \gamma_k)(f(x_{k}) - f^{*}) + \frac{\varepsilon^2}{4LD^2} + \frac{\varepsilon^2}{8LD^2} \\
& = & (1 - \gamma_k)(f(x_{k}) - f^{*}) + \frac{3\varepsilon^2}{8LD^2} \\
& \le & (1 - \gamma_k)(f(x_{k}) - f^{*}) + \frac{3\varepsilon}{8LD^2}(f(x_{k}) - f^{*}) \\
& = & (1 - \frac{\varepsilon}{2LD^2} + \frac{3\varepsilon}{8LD^2})(f(x_{k}) - f^{*}) \\
& = & (1 - \frac{\varepsilon}{8LD^2})(f(x_{k}) - f^{*}),
\end{eqnarray*}
or 
$$f(x_{k+1}) - f(x_k) \le -\frac{\varepsilon^2}{8LD^2}.$$

\textit{Case II}: If $\gamma_k = 1$, then $2LD^2 \le \varepsilon$. Therefore, 
\begin{eqnarray*}
f(x_{k+1}) & \le & f(x_k) + \gamma_k(f^{*} - f(x_k)) + \frac{\gamma_k\varepsilon}{2} + \frac{L\gamma^2_kD^2}{2} \\
& \le & f(x_k) - \gamma_k\varepsilon + \frac{\gamma_k\varepsilon}{2} + \frac{L\gamma^2_kD^2}{2} \\
& \le & f(x_k) - \gamma_k\varepsilon + \frac{\gamma_k\varepsilon}{2} + \frac{\gamma^2_k\varepsilon}{4} \\
& = & f(x_k) - \frac{\varepsilon}{4}.
\end{eqnarray*}

Combining these two cases we have the improvement
$$f(x_{k+1}) - f(x_k) \le -\beta_1 \varepsilon,$$
where $\beta_1 = \min\{\frac{\varepsilon}{8LD^2}, \frac{1}{4}\}$.
\end{proof}

\subsection{Results for the Away-step Frank-Wolfe Algorithm}

\begin{subdef}
Fix some $\epsilon_g \in (0, (4D)^{-1})$, where $D$ is the diameter of the feasible region $\mathcal{X}$. We say that the vector $g_k$ is a {\em good gradient approximation} at iteration $k$ (with error at most $\varepsilon_g$) if
\[
\|g_k - \nabla f(x_k)\| \le \varepsilon_g(-g^T_k(s_k - v_k)),
\]
where $s_k$ and $v_k$ are, respectively, the vertex obtained by solving the LP and the vertex obtained by enumerating over active vertices in Algorithm~\ref{away_fw_alg} at iteration $k$.
\end{subdef}

From the definition of the direction $d_k$ in Algorithm~\ref{away_fw_alg}, we have the bound
$$2g^T_kd_k \le g_k^T(s_k - x_k) + g_k^T(x_k - v_k) = g_k^T(s_k - v_k) \le 0.$$
A good gradient approximation then satisfies
\[
\|g_k - \nabla f(x_k)\| \le \varepsilon_g(-g^T_k(s_k - v_k)) \le 2\varepsilon_g(-g^T_kd_k).
\]

Under the good gradient approximation assumption that $\varepsilon_g \le (4D)^{-1}$, 
and using the fact that $\|d_k\|$ is bounded by the diameter $D$, the Cauchy-Schwarz inequality gives
$$|(g_k - \nabla f(x_k))^Td_k| \le 2\varepsilon_g(-g_k^Td_k)D,$$
implying that
\begin{equation}\label{grad2}
\nabla f(x_k)^Td_k \le g_k^Td_k(1 - 2\varepsilon_gD).
\end{equation}

Now, $g_k^Td_k \le 0$, so \eqref{grad2} with $\varepsilon_g \le (4D)^{-1}$ implies that $\nabla f(x_k)^Td_k \le 0$. Let $x^*$ denote the optimal solution, with objective function value $f^{*} = f(x^{*})$. Then, by the definition of $d_k$, the Cauchy-Schwarz Inequality, and convexity of $f$,
\begin{eqnarray*}
g_k^Td_k & \le & g_k^T(x^{*} - x_k)\\
& = & \nabla f(x_k)^T(x^{*} - x_k) + (g_k - \nabla f(x_k))^T(x^{*} - x_k) \\
& \le & \nabla f(x_k)^T(x^{*} - x_k) + 2\varepsilon_g(-g_k^Td_k)\|x^{*} - x_k\| \\
& \le & f^{*} - f(x_k) + 2\varepsilon_gD(-g_k^Td_k),
\end{eqnarray*}
which gives
\begin{equation}\label{grad3}
g_k^Td_k \le \frac{f^{*} - f(x_k)}{1 + 2\varepsilon_gD}. 
\end{equation}
\noindent
From the step size policy, $$\gamma_k = \min\left\{\gamma^{max}_k, -\frac{g_k^Td_k}{L\|d_k\|^2}\right\},$$
we see that $\gamma_k L\|d_k\|^2 \le -g_k^Td_k$. We are now ready to state and prove two results on the decrease in function value, corresponding to the two cases in the step size policy.

\begin{sublemma} Suppose that we have a good gradient approximation at iteration $k$ and the step size is at its maximum, i.e. $\gamma_k = \gamma^{max}_k$. Then,
$$
f(x_{k+1}) - f^{*} \le \left(1 - \gamma_k\frac{1/2 - 2\varepsilon_gD}{1 + 2\varepsilon_gD}\right)(f(x_k) - f^{*} ).
$$
\label{dec_max_step}
\end{sublemma}
\begin{proof}
By $L-$Lipschitz continuity and \eqref{grad2},
\begin{eqnarray*}
f(x_{k+1}) & \le & f(x_k) + \gamma_k\nabla f(x_k)^Td_k + \frac{\gamma_k^2}{2}L\|d_k\|^2 \\ 
& \le & f(x_k) + \gamma_k \nabla f(x_k)^Td_k - \frac{g_k^Td_k}{2}\gamma_k \,\,\, (\text{by definition of $\gamma_k$}) \\
& \le & f(x_k) + \gamma_k g_k^Td_k\left(1 - 2\varepsilon_gD - 1/2\right) \,\,\, (\text{from \eqref{grad2}}) \\
& = &  f(x_k) + \gamma_k g_k^Td_k(1/2 - 2\varepsilon_gD ). 
\end{eqnarray*}
Using \eqref{grad3}, we have that 
$$f(x_{k+1}) \le f(x_k) + \gamma_k g_k^Td_k(1/2 - 2\varepsilon_gD ) \le f(x_k)  + \gamma_k (f^{*} - f(x_k))\frac{1/2 - 2\varepsilon_gD}{1 + 2\varepsilon_gD}.$$
Subtracting $f^*$ from both sides, we obtain
$$f(x_{k+1}) - f^{*} \le \left(1 - \gamma_k\frac{1/2 - 2\varepsilon_gD}{1 + 2\varepsilon_gD}\right)(f(x_k) - f^{*} ).$$
\end{proof}
 
\begin{sublemma}\label{lem:goodstep1}
Suppose that we have a good gradient approximation at iteration $k$ and the step size is $\gamma_k = -g_k^Td_k / (L\|d_k\|^2)$. Then 
\begin{equation}\label{fn_dec}
f(x_{k+1})  \le f(x_k) - \frac{(g_k^T(s_k - v_k))^2}{4L\|d_k\|^2} \left(\frac{1}{2} - 2\varepsilon_gD \right).
\end{equation}
\end{sublemma}

\begin{proof}
Using \eqref{grad2}, the expression for $\gamma_k$, and the $L-$Lipschitz property,
\begin{eqnarray*}
f(x_{k+1}) & \le & f(x_k) + \gamma_k\nabla f(x_k)^Td_k + \frac{\gamma_k^2}{2}L\|d_k\|^2 \\
& \le & f(x_k) + \gamma_k g_k^Td_k (1 - 2\varepsilon_gD) + \frac{\gamma_k^2}{2}L\|d_k\|^2 \,\,\, (\text{from \eqref{grad2} }) \\
& = & f(x_k) - \frac{(g_k^Td_k)^2}{L\|d_k\|^2}(1 - 2\varepsilon_gD) + \frac{(g_k^Td_k)^2}{2L\|d_k\|^2} \,\,\, (\text{by $\gamma_k$ substitution}) \\
& = & f(x_k) - \frac{(g_k^Td_k)^2}{L\|d_k\|^2}\left(1 - 2\varepsilon_gD - \frac{1}{2} \right).
\end{eqnarray*}
\noindent
Using the fact that $g_k^Td_k \le \frac{1}{2} g_k^T(s_k - v_k) \le 0$, we then have 
$$f(x_{k+1}) \le f(x_k) - \frac{(g_k^T(s_k - v_k))^2}{4L\|d_k\|^2}\left(\frac{1}{2} - 2\varepsilon_gD \right).$$
\end{proof}
Next, we apply some definitions 
and results for polytopes from \cite{beckshtern2017linearly} to the feasible region $\mathcal{X}$. 

\begin{subdef}
Let $I(x)$ be the index set of active constraints of a vector $x$ in the polytope $\mathcal{X} = \{x: Ax \le b\}$, i.e. 
$$I(x) = \{i : A_ix = b_i\},$$
where $A_i$ is the $i$th row of the matrix $A$.\end{subdef}

\begin{subdef}
Similarly, given a set $U \subseteq \mathcal{V}$ of vertices of the polytope $\mathcal{X}$, the set of active constraints for all points in $U$ is 
$$I(U) = \{i : A_iv = b_i \,\, \forall v \in U\} = \bigcap_{v \in U}I(v).$$
\end{subdef}

\begin{sublemma}
Let $x \in \mathcal{X}$ and $U \subseteq \mathcal{V}$ satisfy $x = \sum_{v \in U} \alpha_v v$, where all vertices in $U$ have a positive coefficient $\alpha_v$, and $\sum_{v \in U} \alpha_v = 1$. Then $I(x) = I(U)$.
\label{index_set}
\end{sublemma}

The following technical lemma is from \cite{beckshtern2017linearly} which is based on a polytope's pyramidal width constant. This value is unknown and difficult to compute, even in the case of the unit simplex, but its existence is a boon to showing convergence. For any matrix $B$ and any index set $I$, define $B_I$ to be the submatrix of $B$ consisting of the rows of $B$ whose indices appear in $I$.

\begin{sublemma}
Let $U \subseteq \mathcal{V}$ and $w \in \mathbb{R}^n$ be given. If there exists $z \in \mathbb{R}^n$ such that $A_{I(U)}z \le 0$ and $w^Tz > 0$, then 
$$\max_{p \in \mathcal{V}, u \in U} w^T(p-u) \ge \frac{\Omega_\mathcal{X}}{|U|}\frac{w^Tz}{\|z\|},$$
where $\Omega_{\mathcal{X}} = \zeta / \varphi$
for
\begin{eqnarray*}
\zeta & = & \min_{v \in\mathcal{V}, i: b_i > A_iv} (b_i - A_iv) \text{ and} \\
\varphi & = & \max_{\mathcal{I} \setminus I(\mathcal{V})} \|A_i\|.
\end{eqnarray*}
Here, $\mathcal{I}$ is the set of indices of all constraints. 
\label{pdirw}
\end{sublemma}

Lemma~\ref{lem:goodstep1} gives a bound on expected function improvement, but the bound relies on the alignment between the estimated gradient and the vertices identified in the away-step method. We are now in a position to give a different bound that is central in deriving a rate of convergence in the next section, by relating the bound of Lemma~\ref{lem:goodstep1} to the pyramidal width of the polytope using Lemma~\ref{pdirw}. Recall that $N$ is the number of vertices of the polytope $\mathcal{X}$, $L$ is the Lipschitz constant of the gradient and $\mu$ is the strong convexity parameter.

\begin{sublemma}\label{lem:good_step_2}
Suppose that we have a good gradient approximation at iteration $k$ and the step size is $\gamma_k = -g_k^Td_k / (L\|d_k\|^2)$. Then
$$
f(x_{k+1}) - f^{*}  \le \left(1 - \left(\frac{\Omega_{\mathcal{X}}}{N}\right)^2 \frac{\mu(1/2 - 2\varepsilon_gD)}{8LD^2(2\varepsilon_gD + 1)^2} \right)(f(x_k) - f^{*}).  
$$
\end{sublemma}

\begin{proof}
By convexity,
$$f^{*} \ge f(x_k) + \nabla f(x_k)^T(x^{*} - x_k).$$
The Cauchy-Schwarz Inequality gives 
$$(g_k - \nabla f(x_k))^T(x^{*} - x_k) \le \|g_k - \nabla f(x_k)\|\,\,\|x^{*} - x_k\| \le 2\varepsilon_g(-g_k^Td_k)\|x^{*}-x_k\|.$$

Therefore, 
$$f^{*} \ge f(x_k) + \nabla f(x_k)^T(x^{*} - x_k) \ge f(x_k) + g_k^T(x^{*} - x_k) +  2\varepsilon_g(g_k^Td_k)\|x^{*}-x_k\|.$$
Rearranging gives
$$-g_k^T((x^{*} - x_k) + 2\varepsilon_gd_k\|x^{*}-x_k\|) \ge f(x_k) - f^{*} > 0.$$

Let $z = (x^{*} - x_k) + 2\varepsilon_gd_k\|x^{*}-x_k\|$ and $w = -g_k$. We will prove that $z, w$ satisfy the conditions of Lemma \ref{pdirw}. To that end, $x_k$ is the solution at iteration $k$, and let $U$ be the set of active vertices of $x_k$, i.e. $U = \{u \in \mathcal{V} : \alpha^{(k)}_u > 0\}$. By the definition of the index set and Lemma \ref{index_set}, we have that 
$$A_{I(U)}x_k = A_{I(x_k)}x_k = b_{I(x_k)} = b_{I(U)}.$$
\noindent
Moreover, by the definition of the feasible set,
$A_{I(U)}x^{*} \le b_{I(U)}$
which gives $A_{I(U)}(x^{*} - x_k) \le b_{I(U)} - b_{I(U)} = 0$. Next, consider $A_{I(U)}d_k$.

\textit{Case I: $d_k = s_k - x_k$.} We know that $A_{I(U)}x_k = A_{I(x_k)}x_k = b_{I(x_k)} = b_{I(U)}$ and $A_{I(U)} s_k \le b_{I(U)}$ by feasibility. Then, $A_{I(U)}(s_k - x_k) \le 0$.

\textit{Case II: $d_k = x_k - v_k$.} Because  $v_k \in U$, then $A_{I(U)}v_k = b_{I(U)}$, giving $A_{I(U)}(x_k - v_k) = b_{I(U)} - b_{I(U)} = 0$.

In both cases, we have $A_{I(U)}d_k \le 0$, and so $A_{I(U)}z \le 0$.

Now, by Lemma~\ref{pdirw} using $w = -g_k$ and $z = 2\varepsilon_gd_k\|x^{*} - x_k\| + (x^{*} - x_k)$,
$$\max_{p \in \mathcal{V}, u \in U} -g_k^T(p - u) \ge \frac{\Omega_{\mathcal{X}}}{|U|}\frac{-g_k^Tz}{\|z\|} \ge \frac{\Omega_{\mathcal{X}}}{|U|}\frac{f(x_k) - f^{*}}{\|z\|}.$$

By the definitions of $v_k$ and $s_k$, 
$$g_k^T(v_k - s_k) = \max_{p \in \mathcal{V}, u \in U} -g_k^T(p - u).$$
\noindent
Hence, using the fact that $\|u+v\|^2 \le (\|u\| + \|v\|)^2$ for any vectors $u,v$, and recalling that $N$ is the number of vertices in $\mathcal{V}$, we have
\begin{eqnarray*}
(g_k^T(v_k - s_k))^2 & \ge & \left(\frac{\Omega_{\mathcal{X}}}{N}\right)^2\frac{(f(x_k) - f^{*})^2}{\| 2\varepsilon_gd_k\|x^{*} - x_k\| + (x^{*} - x_k) \|^2} \\ 
& \ge &  \left(\frac{\Omega_{\mathcal{X}}}{N}\right)^2 \frac{(f(x_k) - f^{*})^2}{(\|2\varepsilon_g d_k\|\,\|x^{*}-x_k\| + \|x^{*} - x_k\|)^2} \\
& \ge &  \left(\frac{\Omega_{\mathcal{X}}}{N}\right)^2 \frac{(f(x_k) - f^{*})^2}{(2\varepsilon_g \|d_k\| + 1)^2 \|x^{*} - x_k\|^2} \\
& \ge &  \left(\frac{\Omega_{\mathcal{X}}}{N}\right)^2 \frac{(f(x_k) - f^{*})^2}{(2\varepsilon_gD + 1)^2\|x^{*}-x_k\|^2}.
\end{eqnarray*}

By strong convexity with parameter $\mu$ and the fact that $x^*$ is optimal,
$$
f(x_k) \ge f^{*} + \nabla f(x^{*})^T(x_k - x^{*}) + \frac{\mu}{2}\|x^{*} - x_k\|^2 \ge f^{*} + \frac{\mu}{2}\|x^{*} - x_k\|^2,
$$
giving
$$\frac{f(x_k) - f^{*}}{\|x^{*} - x_k\|^2} \ge \frac{\mu}{2}.$$
Combining this with the earlier inequality we obtain
\begin{equation}
(g_k^T(v_k - s_k))^2  \ge \left(\frac{\Omega_{\mathcal{X}}}{N}\right)^2 \frac{\mu(f(x_k) - f^{*})}{2(2\varepsilon_gD + 1)^2}.
\label{bound_g_vs}
\end{equation}

From \eqref{fn_dec}, we have 
\begin{eqnarray*}
f(x_{k+1}) & \le & f(x_k) - \frac{(g_k^T(s_k - v_k))^2}{4L\|d_k\|^2} \left(\frac{1}{2} - 2\varepsilon_gD \right) \\
& \le & f(x_k) - \left(\frac{\Omega_{\mathcal{X}}}{N}\right)^2 \frac{\mu(f(x_k) - f^{*})}{8LD^2} \frac{1/2 - 2\varepsilon_gD}{(2\varepsilon_gD + 1)^2}.
\end{eqnarray*}
\noindent
Subtracting $f^*$ from both sides, we get the required result that
$$f(x_{k+1}) - f^{*} \le \left(1 - \left(\frac{\Omega_{\mathcal{X}}}{N}\right)^2 \frac{\mu(1/2 - 2\varepsilon_gD)}{8LD^2(2\varepsilon_gD + 1)^2} \right)(f(x_k) - f^{*}).$$
\end{proof}

We can summarize the results from this section on the away-step Frank Wolfe algorithm with the following result on the function decrease on any iteration. We say that there is a vertex drop in iteration $k$ if, in iteration $k$, we take a maximal away step, so that the vertex $v_k$ is dropped from the active set.

\begin{sublemma}
Suppose that we have a good gradient approximation at iteration $k$. If there is no vertex drop, then 
$$f(x_{k+1}) - f^{*} \le (1 - \beta_2)(f(x_k) - f^{*}),$$
where
\[
\beta_2 = \min\left\{\frac{1/2 - 2\varepsilon_gD}{1 + 2\varepsilon_gD} ,\left(\frac{\Omega_{\mathcal{X}}}{N}\right)^2 \frac{\mu(1/2 - 2\varepsilon_gD)}{8LD^2(2\varepsilon_gD + 1)^2}\right\}.
\]
If there is a vertex drop, then $f(x_{k+1}) - f^{*} \le f(x_k) - f^{*}$.
\label{main_dec_lemma}
\end{sublemma}

\begin{proof}
When there is no vertex drop then we have either $\gamma_k = 1$ for the case when we move to a vertex, in which case Lemma~\ref{dec_max_step} gives the result, or $\gamma_k = -g_k^Td_k/(L\|d_k\|^2) < 1$, in which case Lemma \ref{lem:good_step_2} gives the result. If there is a vertex drop (on an away-step), then the step size achieves its maximum, so we apply Lemma \ref{dec_max_step}, noting that  $\gamma_k(1/2 - 2\varepsilon_gD)/(1 + 2\varepsilon_gD) \ge 0$, to conclude that $f(x_{k+1}) - f^{*} \le f(x_k) - f^{*}$.
\end{proof}

\section{Iteration Complexity}
\label{sec:iteration}

We now develop an iteration complexity analysis through the stochastic process that results from using (random) gradient estimates in the Frank-Wolfe algorithms.
When using gradient estimates, $((X_k, G_k, S_k, V_k, \mathcal{U}_k): k \ge 0)$ is a stochastic process where, on iteration $k$, $X_k$ is the current solution, $G_k$ is the gradient estimate, $S_k$ is the vertex obtained by solving the LP, $V_k$ is the candidate ``bad" vertex used for away steps, and $\mathcal{U}_k$ is the set of active vertices, all at $X_k$. We denote realizations of these random quantities by $(x_k, g_k, s_k, v_k, U_k)$.

For $k \ge 0$, let $\mathcal{F}_k$ denote the $\sigma-$algebra generated from $\{(G_i, S_i, V_i, \mathcal{U}_i)\}^{k-1}_{i=0}$ and $\{X_i\}^k_{i=0}$. Our algorithms are assumed to run only up to the stopping time $T_\varepsilon = \inf\{k \ge 1: f(X_k) - f(x^*) \le \varepsilon \}$; for $k > T_\varepsilon$ we (somewhat arbitrarily) define $X_k = X_{T_\varepsilon}$ and similarly for the other random quantities.

A good gradient estimate arises when the gradient estimate is sufficiently close to the true gradient, but the definition of ``sufficient closeness'' must be customized to the problem at hand \cite{blanchetscheinberg2019trust,paqscheinberg2020ls,menickellyscheinberg2018stoch}. In preparation for our definition, let $\varepsilon_g \in (0, (4D)^{-1})$ be a fixed constant, where $D$ is the diameter of the feasible region $\mathcal{X}$.

\begin{subdef} 
We say, at iteration $k$, that $G_k$ is a \textit{good gradient estimate} for the standard Frank-Wolfe algorithm with at least probability $p_g$ if the event
$$E^{(S)}_k = \{\|G_k - \nabla f(X_k)\| \le \varepsilon/4D\}$$
satisfies the condition $\Prob(E_k^{(S)} | \mathcal{F}_{k}) \ge p_g$  on the event $\{k < T_\varepsilon\}$.
\label{good_grad_def_standard}
\end{subdef}

\begin{subdef} 
We say, at iteration $k$, that $G_k$ is a \textit{good gradient estimate} for the Away-step Frank-Wolfe algorithm with at least probability $p_g$ if the event
$$E^{(A)}_k = \{\|G_k - \nabla f(X_k)\| \le \varepsilon_g G^T_k(V_k - S_k)\}$$
satisfies the condition $\Prob(E_k^{(A)} | \mathcal{F}_{k}) \ge p_g$  on the event $\{k < T_\varepsilon\}$.
\label{good_grad_def_away}
\end{subdef}

Consider the stochastic processes $\Phi^{(S)}, \Phi^{(A)}$ obtained from random variables generated from the standard Frank-Wolfe algorithm, and the Away-step Frank-Wolfe algorithm respectively, where
$$\Phi^{(S)}_k = \exp(f(X_{k}) - f^{*}) ,$$
and
$$\Phi^{(A)}_k = \exp(\nu(f(X_{k}) - f^{*}) + (1 - \nu)N_k).$$
(Recall that $N_k$ is the number of vertices in the active vertex set $\mathcal{U}_k$.) The stochastic processes $\Phi^{(S)}, \Phi^{(A)}$ are bounded because the polytope $\mathcal{X}$ is bounded, implying that $f$ attains both its minimum and maximum on $\mathcal{X}$, and also because $\mathcal{X}$ has a finite number of vertices. Define $M = \max\{\max_{x \in \mathcal{X}}|f(x)|,1\} < \infty$. 

\begin{thm}
There exists $p_g \in (0, 1)$, $\nu \in (0, 1)$ and $\delta_S, \delta_A \in (0, 1)$ such that if we get a good gradient estimate at iteration $k$ with probability at least $p_g$, then
\[
I(T_{\varepsilon} > k)\E\left[\frac{\Phi^{(q)}_{k+1}}{\Phi^{(q)}_k} \biggr| \mathcal{F}_{k}\right] \le e^{-\delta_q}.
\]
for $q = S, A$. 
\label{main_thm}
\end{thm}

\begin{proof}
\textit{Part I: }The standard Frank-Wolfe case, i.e., $q = S$.

From Lemma~\ref{standard_fw_improvement}, for $k < T_{\varepsilon}$,
\[
\frac{\Phi^{(S)}_{k+1}}{\Phi^{(S)}_k}  \le
\begin{cases}
 \exp(-\beta_1\varepsilon) & \text{if we have a good gradient estimate, and} \\
 \exp(2M) & \text{if we have a bad gradient estimate.}
\end{cases}
\]

Let $P_k$ be the conditional probability of a good gradient estimate, conditional on $\mathcal{F}_k$. Then, taking conditional expectations, we get
\begin{align*} 
I(T_{\varepsilon} > k)\E\left[\frac{\Phi^{(S)}_{k+1}}{\Phi^{(S)}_k}\biggr|\mathcal{F}_k\right]
&\le I(T_{\varepsilon} > k) [\exp(2M) - P_k(\exp(2M) - \exp(-\beta_1\varepsilon ))] \notag\\
& \le I(T_{\varepsilon} > k) [\exp(2M) - p_g(\exp(2M) - \exp(-\beta_1\varepsilon))],
\end{align*}
where $p_g$ is the lower bound of the probability of a good gradient approximation. Now, fix $\delta_S$ such that $0 < \delta_S < \beta_1\varepsilon $, which then allows us to select $p_g$ so that
\begin{equation}\label{eq:pglowerFW}
1 > p_g \ge \frac{\exp(2M) - \exp(-\delta_S)}{\exp(2M) - \exp(-\beta_1\varepsilon)} > 0.
\end{equation}
Hence,
$$I(T_{\varepsilon} > k)\E\left[\frac{\Phi^{(S)}_{k+1}}{\Phi^{(S)}_k}\biggr|\mathcal{F}_k\right] \le e^{-\delta_S}.$$

\textit{Part II: } The Away-step Frank-Wolfe case, i.e., $q = A$.

We condition on $\mathcal{F}_k$, so that $\{G_i, S_i, V_i\}^{k-1}_{i = 0}$ and $\{X_i\}^k_{i=0}$ are known. Thus, the active set $\mathcal{U}_k$ and the number of active vertices $N_k$ are also known, but $G_k, S_k$ and $V_k$ are not yet available. Suppose we are yet to stop so that $T_\varepsilon > k$ (which is also $\mathcal{F}_k$ measurable).

Consider the number of active vertices $N_k$. Then
\[
N_{k+1} - N_k = 
\begin{cases}
1 \,\, & \text{if no vertex drop and a new active vertex,} \\
1 - N_k & \text{if no vertex drop and move to a vertex, and}\\
-1\,\, & \text{if vertex drop.}
\end{cases}
\]

By Lemma \ref{main_dec_lemma},
\[
\frac{\Phi^{(A)}_{k+1}}{\Phi^{(A)}_k}  \le
\begin{cases}
 \exp(\nu(-\beta_2\varepsilon) + (1-\nu)({N_{k+1} - N_k})) & \text{if good gradient estimate and no vertex drop,} \\
 \exp(-(1 - \nu)) & \text{if good gradient estimate and vertex drop, and} \\
 \exp(2\nu M + (1 - \nu)) & \text{if bad gradient estimate.}
\end{cases}
\]

Fix $\nu \in (0,1)$ so that $(1 + \beta_2\varepsilon)^{-1} < \nu < 1$, and again let $P_k$ be the conditional probability of a good gradient estimate, conditional on $\mathcal{F}_k$. Then, taking conditional expectations, on the event $T_{\varepsilon} > k$,
\begin{align} 
\lefteqn{\E\left[\frac{\Phi^{(A)}_{k+1}}{\Phi^{(A)}_k}\biggr|\mathcal{F}_k\right]} \notag\\
& \le P_k \max\{\exp(-\nu\beta_2\varepsilon + 1 - \nu), \exp(-(1 - \nu))\} + (1-P_k)\exp(2M\nu + 1 - \nu) \notag\\
&= \exp(2M\nu + 1 - \nu) - P_k(\exp(2M\nu + 1 - \nu) - \max\{\exp(-\nu\beta_2\varepsilon + 1 - \nu), \exp(-(1 - \nu))\}) \notag\\
& \le \exp(2M\nu + 1 - \nu) - p_g(\exp(2M\nu + 1 - \nu) - \max\{\exp(-\nu\beta_2\varepsilon + 1 - \nu), \exp(-(1 - \nu))\}).\label{expected_dec}
\end{align}

Now, fix $\delta_A$ such that $0 < \delta_A < \min\{\nu\beta_2\varepsilon - 1 + \nu, 1 - \nu\}$, which then allows us to select $p_g$ so that
\begin{equation}\label{eq:pglowerASFW}
1 > p_g \ge \frac{\exp(2M\nu + 1 - \nu) - \exp(-\delta_A)}{\exp(2M\nu + 1 - \nu) - \max\{\exp(-\nu\beta_2\varepsilon + 1 - \nu), \exp(-(1 - \nu))\}} > 0.
\end{equation}

Hence, by \eqref{expected_dec},
$$I(T_{\varepsilon} > k)\E\left[\frac{\Phi^{(A)}_{k+1}}{\Phi^{(A)}_k}\biggr|\mathcal{F}_k\right] \le e^{-\delta_A}.$$
\end{proof}

The proof of Theorem~\ref{main_thm} shows that there are many values of $\delta_S, \delta_A$ that satisfy the needed conditions. For example, for the standard Frank-Wolfe algorithm, we can pick $\delta_S = \beta_1\varepsilon / 2$, and for the away-step Frank-Wolfe algorithm we can pick 
\begin{equation}\label{eq:ASFWParams}
\nu = \frac{1}{1 + \beta_2\varepsilon/2}, \delta_A = \frac{\beta_2\varepsilon/2}{2 + \beta_2\varepsilon}.
\end{equation}
We can then obtain the following conclusion on iteration complexity.

\begin{thm} \label{thm:iteration_complexity}
Suppose the initial solution $X_0 = x_0$ is deterministic with $f(x_0) - f^* > \varepsilon$. Let $T_{\varepsilon}$ be the stopping time defined as $T_{\varepsilon} = \inf\{k : f(X_k) - f^{*} \le \varepsilon\}$. The expected stopping times of the standard and away-step Frank-Wolfe algorithms are $\mathcal{O}(\varepsilon^{-2}(f(x_0) - f^*))$ and $\mathcal{O}(\varepsilon^{-1}(f(x_0) - f^*))$ respectively. In addition, for both algorithms the stopping time $T_{\varepsilon}$ has a finite moment generating function in a neighborhood of 0. 
\end{thm}

\begin{proof}
Consider a stochastic process $\Phi = (\Phi_k: k \ge 0)$, adapted to $(\mathcal{F}_k: k \ge 0)$ that satisfies the expected ratio decrease as in Theorem~\ref{main_thm} and define a stochastic process $(M_k: k \ge 0)$ where $M_k = e^{k\delta}\Phi_k I(k \le T_{\varepsilon})$. Then, on the event $k < T_{\varepsilon},$ for $\delta = \delta_S$ or $\delta_A$ for the standard or away-step stochastic process respectively,
\begin{eqnarray*}
\mathbb{E}[M_{k+1} | \mathcal{F}_k] & = & e^{(k+1)\delta}\mathbb{E}[\Phi_{k+1}I(k+1 \le T_{\varepsilon}) | \mathcal{F}_k] \\
& = & e^{(k+1)\delta}\mathbb{E}[\Phi_{k+1}I(k < T_{\varepsilon}) | \mathcal{F}_k] \\
& = & e^{(k+1)\delta}I(k < T_{\varepsilon})\mathbb{E}[\Phi_{k+1} | \mathcal{F}_k] \\
& \le & e^{k\delta}I(k < T_{\varepsilon})\Phi_k = M_k.
\end{eqnarray*}

On the event $T_{\varepsilon} \le k$, $M_{k+1} = 0 \le M_k.$ So $(M_k : k \ge 0)$ is a non-negative super-martingale. By the optional stopping theorem, and since $X_0 = x_0$ is deterministic with $f(x_0) > f^* + \varepsilon$,
$$\mathbb{E}[M_{T_{\varepsilon}}] \le \mathbb{E}[M_0] = \Phi_0.$$
So, $\E[e^{\delta T_\varepsilon}\Phi_{T_{\varepsilon}}] \le \Phi_0$, and, since $\Phi_k\ge 1$ for all non-negative integers $k$, $\E[e^{\delta T_\varepsilon}] \le \Phi_0 < \infty$, establishing the finiteness of the moment generating function in a neighborhood of 0.

By concavity of the logarithm function, 
$$\E[\delta T_{\varepsilon}] \le \log \E[e^{\delta T_{\varepsilon}}] \le \log \Phi_0.$$
Therefore, for the standard Frank-Wolfe algorithm, using Lemma~\ref{standard_fw_improvement},
$$\mathbb{E}[T_{\varepsilon}] \le \frac{\log \Phi^{(S)}_0}{\delta_S} = \frac{2\log \Phi^{(S)}_0}{\beta_1\varepsilon} = 2(f(x_0) - f^*)\max\left\{\frac{8LD^2}{\varepsilon^2}, \frac{4}{\varepsilon}\right\} \sim \mathcal{O}\left(\frac{f(x_0) - f^*}{\varepsilon^2}\right).$$

Similarly, for the away-step Frank-Wolfe algorithm, $\log \Phi^{(A)}_0 = \nu(f(x_0) - f(x^*)) + (1-\nu)|U_0|$ where $|U_0|$ is the initial number of active vertices. For our choice of $\nu$, $1-\nu = \mathcal{O}(\varepsilon)$, so that
$$\mathbb{E}[T_{\varepsilon}] \le \frac{\log\Phi^{(A)}_0}{\delta_A} = (\log\Phi^{(A)}_0)\left(2 + \frac{4}{\beta_2\varepsilon}\right) \sim \mathcal{O}\left(\frac{f(x_0) - f^* + |U_0|}{\varepsilon} \right) = \mathcal{O}\left(\frac{f(x_0) - f^*}{\varepsilon}\right).$$
Here, $\beta_2$ does not depend on $\varepsilon$. 
\end{proof}

\section{Sample Size Complexity}
\label{sec:samplesize}

We establish results on sample-size complexity under two different assumptions, the first of which is more general than the second but leads to a weaker complexity result.

\subsection{The Bounded Second Moment Case}

Recall that, for each feasible $x$, we assume the availability of i.i.d. random vectors $G_1(x), G_2(x), \ldots$ that can be used to estimate the (true) gradient $\nabla f(x)$.

\begin{assp}\label{ass:BoundedSecondMoment}
For all feasible $x$, the random vector $G_1(x)$ is unbiased and has bounded (in $x$) second moment, i.e., $\E G_1(x) = \nabla f(x)$, and there exists a constant $V_g < \infty$ such that, for all $x \in \mathcal{X}$,
$$E[\|G_1(x) - \nabla f(x)\|^2] \le V_g.$$
\end{assp}

The gradient estimator at any point, $x$, is simply the sample average, $g(x) = n^{-1} \sum_{i=1}^n G_i(x)$. We adopt a constant sample size, $n$, at any iterate. Let $\|\cdot\|$ denote the 2-norm.

\begin{sublemma}
Under Assumption~\ref{ass:BoundedSecondMoment}, for any fixed feasible $x$ and any $s > 0$,
$$\Prob(\|g(x) - \nabla f(x)\| > s) \le \frac{V_g}{n s^2}.$$
\label{sample1}
\end{sublemma}

\begin{proof}
Let $\nabla f(x, i), g(x, i)$ and $G(x, i)$ denote the $i$th components of $\nabla f(x), g(x)$ and $G(x)$ respectively. By Chebyshev's inequality,
\begin{align*}
\Pr(\|g(x) - \nabla f(x)\| > s) & \le s^{-2}\E[\|g(x) - \nabla f(x)\|^2]\\
& =  s^{-2}\sum^d_{i=1} \E[(g(x, i) - \nabla f(x, i))^2] \\
& =  \frac{\sum^d_{i=1} \E[(G(x, i) - \nabla f(x, i))^2]}{n s^2} \\
& =  \frac{\E[\|G(x) - \nabla f(x)\|^2]}{n s^2}\\
& \le  \frac{V_g}{n s^2}.
\end{align*}
\end{proof}
Theorem~\ref{sample_size_thm} establishes sample size requirements for the standard and away-step Frank-Wolfe algorithms. The away-step Frank-Wolfe algorithm permits a larger order of error in the gradient estimate than the standard Frank-Wolfe algorithm, and therefore requires a smaller sample size.

\begin{thm}
Under Assumption~\ref{ass:BoundedSecondMoment}, the sample sizes in each iteration required to achieve a \textit{good gradient estimate} for the standard and Away-step Frank-Wolfe algorithms are $\mathcal{O}(\varepsilon^{-4})$ and $\mathcal{O}(\varepsilon^{-2})$, respectively.
\label{sample_size_thm}
\end{thm}

\begin{proof}

By Lemma~\ref{sample1} and Definition~\ref{good_grad_def_standard}, take $s = \varepsilon/(4D)$, and then it suffices to choose the sample size $n$ so that
\begin{equation}\label{eq:nlower-FW}
n \ge \frac{16V_gD^2}{\varepsilon^2(1 - p_g)}.
\end{equation}
Consider the lower bound on $1-p_g$ obtained from \eqref{eq:pglowerFW} with $\delta_S = \beta_1 \epsilon/2$. For any $M \ge 1$, $e^{2M+1} - 1 \le e^{2M}(e^2/2) = e^{2M + 2}/2$, so by the Taylor expansion of $e^{\beta_1\varepsilon/2}$,
\begin{equation}\label{eq:LowerBdOn1-pg}
1 - p_g = \frac{\exp(-\beta_1\varepsilon/2) - \exp(-\beta_1\varepsilon)}{\exp(2M) - \exp(-\beta_1\varepsilon)} = \frac{\exp(\beta_1\varepsilon/2) - 1}{\exp(2M + \beta_1\varepsilon) - 1} \ge \frac{\beta_1\varepsilon/2}{e^{2M + 2}/2}.
\end{equation}
Hence, the necessary sample size $n$ in \eqref{eq:nlower-FW} is
$$\frac{16V_gD^2}{\varepsilon^2(1 - p_g)} \le \frac{16V_gD^2e^{2M+2}}{\beta_1\varepsilon^3} = \frac{128LV_gD^4e^{2M+2}}{\varepsilon^4} = \mathcal{O}(\varepsilon^{-4}). $$
This completes the proof for the standard Frank-Wolfe algorithm.

Next, consider the away-step Frank-Wolfe algorithm. In view of Definition~\ref{good_grad_def_away}, we want $s$ in Lemma~\ref{sample1} to be a lower bound on $|\varepsilon_g G^T_k(V_k - S_k)|$. From \eqref{bound_g_vs}, it suffices to take 
\begin{equation}\label{eq:ASFW_s}
s^2=\frac{\Omega_{\mathcal{X}}^2}{N^2} \frac{\mu}{2(2\varepsilon_gD + 1)^2}\varepsilon.
\end{equation}
From Lemma ~\ref{sample1}, we obtain a good gradient estimate with sample size
\begin{equation}\label{eq:n_lower_AFW}
n \ge \frac{2V_g(2\varepsilon_gD + 1)^2}{(1-p_g)\varepsilon}\left(\frac{N}{\Omega_{\mathcal{X}}}\right)^2.
\end{equation}
With the choice of $\nu$ and $\delta_A$ in \eqref{eq:ASFWParams}, and using Taylor expansions with sufficiently small $\varepsilon$, 
\begin{equation}\label{eq:AS_low_1-pg}
1 - p_g = \frac{\exp\left(\frac{\beta_2\varepsilon/4}{1 + \beta_2\varepsilon/2}\right) - 1}{\exp\left(\frac{2M+\beta_2\varepsilon}{1+\beta_2\varepsilon/2}\right) - 1} \ge \frac{\exp(\beta_2\varepsilon/6) - 1}{\exp(2M) - 1} \ge \frac{\beta_2\varepsilon/6}{\exp(2M) - 1}\ge e^{-2M-2}\beta_2\varepsilon.
\end{equation}
Then, the necessary sample size in \eqref{eq:n_lower_AFW} is  
$$ \frac{2V_g(2\varepsilon_gD + 1)^2}{(1-p_g)\varepsilon}\left(\frac{N}{\Omega_{\mathcal{X}}}\right)^2 \le \frac{2V_g(2\varepsilon_gD + 1)^2e^{2M+2}}{\beta_2\varepsilon^2}\left(\frac{N}{\Omega_{\mathcal{X}}}\right)^2 = \mathcal{O}(\varepsilon^{-2}).$$
\end{proof}

\subsection{The Sub-Gaussian Case}

Under the bounded variance assumption we were able to establish a rate of convergence of the standard and away-step algorithms. Those rates can be shown to be faster if we make the stronger assumption that the vector $G_1(x)$ that represents a single gradient realization at the feasible point $x$ is sub-Gaussian with a parameter that is uniform in $x$, following \cite{jinetal2019subG, veshynin2012randmat}.

\begin{subdef}\label{def:subGaussian}
A random vector $X \in \mathbb{R}^d$ is sub-Gaussian with parameter $\sigma$ if
$$\mathbb{E}[e^{v^T(X - \mathbb{E}[X])}] \le e^{\frac{\|v\|^2\sigma^2}{2}} \ \ \forall v \in \mathbb{R}^d.$$
\end{subdef}

\begin{assp}
\label{subGassp}
For all feasible $x$, the random vector $G_1(x)$, representing a gradient realization at the point $x$, is unbiased and sub-Gaussian with parameter $\rho/\sqrt{d}$, where $d$ is the dimension of $x$.
\end{assp}

Under Assumption~\ref{subGassp}, $G_1(x) - \nabla f(x)$ is \textit{sub-Gaussian} with parameter $\rho$, uniformly in all feasible points $x$ \cite{jinetal2019subG}. 

\begin{sublemma}
    Let $g(x) = n^{-1} \sum_{i=1}^n G_i(x)$ be a sample average of i.i.d. gradient realizations at $x$. Under Assumption~\ref{subGassp}, there exists a constant $c > 0$ such that, for all feasible $x$,
    $Pr(||g(x) - \nabla f(x)|| \ge s) \le 2de^{-ncs^2}$,
    where $d$ is the dimension of $x$. 
\label{gauss_prob_bound}
\end{sublemma}

\begin{proof}
Fix $s > 0$. Corollary 7 in \cite{jinetal2019subG} establishes that, for some positive constant $\kappa$,
$$Pr\left(\left\|\sum^n_{i=1}(G_i(x) - \nabla f(x))\right\| \le \kappa \sqrt{n\rho^2 \log\frac{2d}{\delta}}\right) \ge 1 - \delta$$
for any $\delta > 0$. Set $\delta = 2d \exp(-ns^2\kappa^{-2}\rho^{-2})$ to obtain the desired result, where $c =\kappa^{-2}\rho^{-2}$. 
\end{proof}

Assumption~\ref{subGassp} enables us to improve the sample size conclusions in Theorem~\ref{sample_size_thm} where we made the weaker assumption of bounded gradient variance.

\begin{thm}\label{th:SampleSizeSubG}
    Under Assumption~\ref{subGassp}, the sample sizes in each iteration required to achieve a \textit{good gradient approximation} for the standard and Away-step Frank-Wolfe algorithms are $\mathcal{O}(\varepsilon^{-2}\log(\varepsilon^{-1}))$ and $\mathcal{O}(\varepsilon^{-1}\log(\varepsilon^{-1}))$, respectively. 
\end{thm}

\begin{proof}
For the standard Frank-Wolfe algorithm, take $s = \varepsilon/(4D)$ in Lemma~\ref{gauss_prob_bound}. Now choose $n$ so large that $2de^{-nc\varepsilon^2/(16D^2)}$ is a lower bound on \eqref{eq:LowerBdOn1-pg}, i.e.
$$n \ge \frac{16D^2}{c\varepsilon^2}(2M + 2 + \log 2d) + \frac{16D^2}{c\varepsilon^2}\log\left(\frac{1}{\beta_1\varepsilon}\right) \sim \mathcal{O}\left(\frac{1}{\varepsilon^2}\log\frac{1}{\varepsilon}\right),$$
and this gives the result for the standard Frank-Wolfe algorithm.

Now, consider the Away-step Frank-Wolfe algorithm. Take $s^2 = c_1\varepsilon$ as in \eqref{eq:ASFW_s} in Lemma~\ref{gauss_prob_bound}. It then suffices to take $n$ so large that $2d\exp(-ncc_1\varepsilon) \le 1 - p_g$. Using the bound \eqref{eq:AS_low_1-pg}, we get 
$$n \ge \frac{2M+2 + \log 2d - \log \beta_2 \varepsilon}{c c_1 \varepsilon} \sim \mathcal{O}\left(\frac{1}{\varepsilon}\log\frac{1}{\varepsilon}\right).$$

\end{proof}

\bibliographystyle{siam}
\bibliography{FWbib}

\end{document}